\documentclass[11pt,reqno,a4paper]{amsart}%
\usepackage{amsfonts,amsmath,times}
\usepackage[scaled=0.92]{helvet}
\usepackage{tikz}

\usepackage{a4wide}
\allowdisplaybreaks

\definecolor{blau}{rgb}{0,0,0.75} 
\usepackage[pdftex]{hyperref}
\hypersetup{colorlinks,linkcolor=blau,citecolor=blue,urlcolor=blau}

\newcommand{\fallfak}[2]{\ensuremath{#1^{\underline{#2}}}}

\newcommand{\N}{\ensuremath{\mathbb{N}}}

\def\P{{\mathbb {P}}}
\def\E{{\mathbb {E}}}

\def\N{{\mathbb {N}}}

\newtheorem{theorem}{Theorem}
\newtheorem{lemma}{Lemma}

\theoremstyle{remark}

\title{On death processes and urn models}

\author[M.~Kuba]{Markus Kuba}
\address{Markus Kuba\\
Institut f{\"u}r Diskrete Mathematik und Geometrie\\
Technische Universit\"at Wien\\
Wiedner Hauptstr. 8-10/104\\
1040 Wien -- Institut f\"ur Angewandte Mathematik und Naturwissenschaften, FH Technikum Wien, H\"ochst\"adtplatz 5, 1200 Wien, Austria} %
\email{kuba@dmg.tuwien.ac.at}

\author[A.~Panholzer]{Alois Panholzer}
\address{Alois Panholzer\\
Institut f{\"u}r Diskrete Mathematik und Geometrie\\
Technische Universit\"at Wien\\
Wiedner Hauptstr. 8-10/104\\
1040 Wien, Austria} \email{Alois.Panholzer@tuwien.ac.at}

\thanks{The second author was supported by the Austrian Science Foundation FWF, grant S9608-N13.}

\date{\today}


\begin{document}

\begin{abstract}
We use death processes and embeddings into continuous time in order to analyze several urn models with a diminishing content. 
In particular we discuss generalizations of the pill's problem, originally introduced by Knuth and McCarthy, 
and generalizations of the well known sampling without replacement urn models, and OK Corral urn models.
 \end{abstract}

\maketitle

\section{Introduction}
\subsection{Diminishing P\'olya-Eggenberger urn models}
In this work we are concerned with so-called P{\'{o}}lya-Eggenberger urn
models, which in the simplest case of two types of colors for the
balls can be described as follows. At the beginning, the urn
contains $n$ white and $m$ black balls. At every step, we choose a
ball at random from the urn, examine its color and put it back into
the urn and then add/remove balls according to its color by the
following rules. If the ball is white, then we put $\alpha$ white and $\beta$
black balls into the urn, while if the ball is black, then $\gamma$ white
balls and $\delta$ black balls are put into the urn. The values $\alpha, \beta, \gamma,
\delta \in \mathbb{Z}$ are fixed integer values and the urn model is
specified by the transition matrix $M = \bigl(\begin{smallmatrix} \alpha &  \beta \\
\gamma & \delta \end{smallmatrix}\bigr)$. Urn models with $r$ ($\ge2$) types of
colors can be described in an analogous way and are specified by an
$r \times r$ transition matrix. Urn models are simple, useful mathematical tools for describing many
evolutionary processes in diverse fields of application such as
analysis of algorithms and data structures, statistics and genetics.
Due to their importance in applications, there is a huge literature
on the stochastic behavior of urn models; see for example
\cite{JohKot1977,KotBal1997,Mah2003}. Recently, a few different approaches
have been proposed, which yield deep and far-reaching results \cite{bai2002,FlaDumPuy2006,FlaGabPek2005,Jan2004,Jan2006,Pou2008}.
Most papers in the literature impose the so-called tenability
condition on the transition matrix, so that the process can be
continued ad infinitum, or no balls of a given color being
completely removed. However, in some applications, examples given
below, there are urn models with a very different nature, which we
will refer to as ``\emph{diminishing urn models}.'' 
Such models have recently received some attention, see for example~\cite{PuyPhD,FlaHui2008,FlaDumPuy2006,Stad1998,BrePro2003,WilMcI1998}. 
For simplicity of presentation, we describe diminishing urn models in the case of balls with two types
of colors, black and white. 

\smallskip

We consider P{\'{o}}lya-Eggenberger urn
models specified by a transition matrix $M = \bigl(\begin{smallmatrix}\alpha &  \beta \\
\gamma & \delta\end{smallmatrix}\bigr)$,
and in addition we also specify a set of absorbing states $\mathcal{A}
\subseteq \mathbb{N} \times \mathbb{N}$. The evolution of the urn takes place
in the state space $\mathcal{S}\subseteq \mathbb{N} \times \mathbb{N}$.
The urn contains $m$ black balls and $n$ white balls at the beginning, with $(m,n)\in\mathcal{S}$, and evolves by successive
draws at discrete instances according to the transition matrix until
an absorbing state in $\mathcal{A}$ is reached, and the process stops. 
Diminishing urn models with more than two type of balls can be
considered similarly. 

\subsection{Plan of this note and notations}
There are numerous examples of diminishing urns and related problems in literature. In the following we present three concrete problems, the pills problem, the sampling without replacement urn, and the OK Corral urn model, and summarize known results. For all three problems presented below, and suitable generalization of them, we will use stochastic processes and an embedding in continuous time, in order to unify and extend 
known results in the literature concerning exact distribution laws, generalizing some results of~\cite{BrePro2003,HwaKuPa2007,PanKu2007+,Stad1998,Kin1999,Kin2002,KinVol2003,PuyPhD,FlaDumPuy2006} .
We will denote with $X\oplus Y$ the sum of independent random variable $X$ and $Y$. Moreover, we use the notations $\N=\{1,2,3,\dots\}$ and $\N_0=\{0,1,2,\dots\}$.

\subsection{The pills problem and generalizations}
Consider the diminishing urn problem with transition matrix
given by $M = \bigl(\begin{smallmatrix} -1 & 0 \\ 1 &
-1\end{smallmatrix}\bigr)$, state space $\mathcal{S}=\N\times\N$, and the absorbing axis $\mathcal{A} =
\{(0,n)\mid n\in\N\}$. An interpretation is as
follows. An urn has two
types of pills in it, which are single-unit and double-unit pills,
respectively. At every step, we pick a pill uniformly at random. If
a single-unit pill is chosen, then we eat it up, and if the pill is
of double unit, we break it into two halves---one half is eaten up
and the other half is now considered of single unit and thrown back
into the urn. The question is then, when starting with $n$
single-unit pills and $m$ double-unit pills, what is the probability
that $k$ single-unit pills remain in the urn when all double-unit
pills are drawn? This problem has been stated by Knuth and McCarthy in~\cite{KnuMcC1991},
where the authors asked for a formula for the expected number of remaining single-unit
pills, when there are no double-unit pills in the urn. The solution
appeared in~\cite{Hes1992}. A more refined study was given by Brennan and Prodinger~\cite{BrePro2003}, where they derive exact formul{\ae} for the
variance and the third moment of the number of remaining single-unit
pills; furthermore, a few generalizations are proposed. The probability generating functions and limit laws for the pills problem and a variant of the problem have been derived in~\cite{HwaKuPa2007} using a generating functions approach. Furthermore, a study of the arising limiting distributions of a general class of related problems has been carried out in~\cite{PanKu2007+}
using a recursive approach basically guessing the structure of the moments, together with an application of the so-called method of moments. 
However, some cases proved to be quite elusive using the techniques of~\cite{PanKu2007+}. Moreover, no simple explicit general formula for the probability mass function of the random variable of interest was obtain 
before. Furthermore, we will discuss in this note weighted generalizations of the pills problem,

\subsection{Sampling without replacement and generalizations}
This classical example, often severing as a toy model, corresponds to the urn with transition matrix $M =
\bigl(\begin{smallmatrix} -1 & 0 \\ 0 & -1\end{smallmatrix}\bigr)$
with absorbing axis: $\mathcal{A} = \{(0,n) \mid n \in\N\}$. In this model, balls are drawn one
after another from an urn containing balls of two different colors
and not replaced. What is the probability that $k$ white balls color remain when black balls have been removed, 
starting with $n$ white and $m$ black balls? This simple urn model has been discussed in~\cite{HwaKuPa2007} using generating functions. 
Moreover, generalizations of the sampling urn model have been discussed in~\cite{PanKu2007+,Me2010}. 

\subsection{The OK Corral urn model} The so-called OK Corral urn serves as a mathematical model of the historical gun fight at the OK Corral.
This problem was introduced by Williams and McIlroy in \cite{WilMcI1998} and studied recently by several authors using
different approaches, leading to very deep and interesting results; see \cite{Stad1998,Kin1999,Kin2002,KinVol2003,PuyPhD,FlaDumPuy2006}. Also the urn corresponding
to the OK corral problem can be viewed as a basic model in the mathematical theory of warfare and conflicts; see \cite{Kin2002,KinVol2003}.

In the diminishing urn setting the OK corral problem corresponds to the urn with transition matrix 
$M =\bigl(\begin{smallmatrix} 0 & -1 \\ -1 & 0\end{smallmatrix}\bigr)$ with two absorbing axes: $\mathcal{A} =
\{(0,n) \mid  n\in\N\} \cup \{(m,0) \mid m\in\N\}$. An interpretation is as
follows. Two groups of gunmen, group A and group B (with $n$ and $m$
gunmen, respectively), face each other. At every discrete time step,
one gunman is chosen uniformly at random who then shoots and kills
exactly one gunman of the other group. The gunfight ends when
one group gets completely ``eliminated". Several questions are of
interest: what is the probability that group A (group B)
survives, and what is the probability that the gunfight ends
with $k$ survivors of group A (group B)? This model was analyzed by Williams and McIlroy~\cite{WilMcI1998}, who obtained an interesting result for the expected value of the number of survivors. Using martingale arguments and the method of moments Kingman~\cite{Kin1999} gave limiting distribution results for the OK Corral urn model
for the total number of survivors. Moreover, Kingman~\cite{Kin2002} obtained further results in a very general setting of Lanchester's theory of warfare. Kingman and Volkov~\cite{KinVol2003} gave a more detailed analysis of the balanced OK Corral urn model using a connection to the famous Friedman urn model; amongst others, they derived an explicit result for the number of survivors and even local limit laws. 
In his Ph.~D.~thesis~\cite{PuyPhD} Puyhaubert extended the results of \cite{Kin1999,KinVol2003} on the balanced OK Corral urn model using analytic combinatoric methods concerning the number of survivors of a certain group. His study is based on the connection to the Friedman urn showed in~\cite{KinVol2003}. He obtained explicit expression for the probability distribution, the moments, and also reobtained (and refined) most of the limiting distribution results reported earlier. Some results of~\cite{PuyPhD} where reported in the work of Flajolet et al.~\cite{FlaDumPuy2006}.
Apparently unknown to the previously stated authors was the earlier work of Stadje~\cite{Stad1998}, who obtained several limiting distribution results for the generalized OK Corral urn, as introduced below, and also for related urn models with more general transition probabilities.  In~\cite{Stad1998} the probability distributions for the most general transition probabilities are ingeniously determined by a complex integral.
The results of Stadje were then discussed in~\cite{Me2010}, and their connection to sampling without replacement, a duality relation, uncovered. However, no transparent probabilistic
derivation of the results of~\cite{Me2010} were given before.

\section{Probabilistic analysis of the pill's problem urn models}
We are interested in a generalized pill's problem with ball replacement matrix given by
\begin{equation}
\label{AofA2012-Pills}
M =\bigl(\begin{smallmatrix} -\alpha & 0 \\ \gamma &-\delta\end{smallmatrix}\bigr),
\end{equation}
where $\alpha,\beta\in\N$, and $\gamma=\alpha\cdot p$, $p\in\N_0$. Let $X_{n,m}$ denote the random variable counting the number of remaining white balls (divided by $\alpha$) when all black balls have been drawn.
The probability generating function $h_{n,m}(v)=\E(v^{X_{n,m}})$ satisfies the recurrence relation
\begin{equation}
\label{AofA2012-PillsRec}
h_{n,m}(v) = \frac{\alpha n}{\alpha n+\delta m}h_{n-1,m}(v) + \frac{\delta m}{\alpha n+\delta m}h_{n+p,m-1}(v),
\end{equation}
with $h_{n,0}=v^n$, $n\ge 1$. We analyze $X_{n,m}$ using a continuous time embedding. We start at time zero with $n$ white balls and $m$ black balls, and use two independent linear processes.
The first one consists of $n$ independent ordinary death processes (white balls) with death rate $\alpha$. 
Let $W_n(t)$ denote the random variable counting the number of living white balls at time $t$, with $W_n(0)=n$. The second one (black balls) consists of $m$ independent modified death processes, with rate $\delta$, 
where each black ball gives at his death birth to $p$ new white balls with death rate $\alpha$, independent of all other balls, and $p\in\N_0$. 
We denote with $B_m(t)$ the random variable counting the number of living black balls at time $t$, with $B_m(0)=m$
Finally, let $C_m(t)$ denote the random variable counting the number of surviving white balls up to time $t$, which are children of black balls.
Let $\tau=\inf_{t>0}\{B_t=0\}$ be the time when the black balls die out. Then
\[
X_{n,m}=W_{n}(\tau)\oplus C_m(\tau),
\]
both random variables are independent, due to the construction. One readily obtains the recurrence relation for the probability generating function $h_{n,m}(v)$
by looking at the time when the first particle dies. It is well known that the index of the variable achieving the minimum out of $r$ independent exponential distributed random variables $X_1,\dots,X_r$ with parameters $\lambda_1,\dots,\lambda_r$,
is given by 
\[
\P\{X_k=\min\{X_1,\dots,X_r\}\}=\int_{0}^{\infty} \Big( \frac{d }{dt} (1-e^{-\lambda_k t}\Big) \prod_{\substack{j=1\\j\neq k}}e^{-\lambda_jt}dt
=\frac{\lambda_k}{\lambda_1+\dots+\lambda_k}.
\]
Hence, we easily obtain that the probability that any of the type one balls dies first is given by 
$\frac{\alpha n}{\alpha n+\delta m}$, and the opposite case happens with probability
$\frac{\delta m}{\alpha n+\delta m}$. Moreover, if $p$ new white balls with death rate $a$ are being born, 
they can be grouped with the already existing white balls, due to the memorylessness 
of exponential distributions. This leads directly to~\eqref{AofA2012-PillsRec}.

\smallskip

Due to the construction of the two processes the random variables $W_n(t)$ and $C_m(t)$
can be decomposed themselves into sums of i.~i.~d. random variables, 
\[
W_n(t)=\bigoplus_{k=1}^{n}X_k(t),\qquad C_m(t)=\bigoplus_{k=1}^{m}Y_k(t)
\]
where $X_k(t)$ denotes the indicator variable of the $k$-th white ball living at time $t$, $1\le k\le n$, 
and $Y_k(t)$ denote the random variable counting the number of surviving white balls up to time $t$, which are children of the $k$-th black ball, $1\le k\le m$.

\smallskip

The probability generating function of a single white ball at time $t$ with rate $\delta$ 
is given by
\[
\E(v^{X_k(t)})=\P\{X_k(t)\le t\} + v\P\{X_k(t)> t\}=1-e^{-\alpha t} + ve^{-\alpha t}.
\]
so the probability generating function of the total number of white balls living at time $t$ is given by 
\[
\E(v^{W_n(t)})=\prod_{k=1}^{n}\E(v^{X_k(t)})=\big(1+(v-1)e^{-\alpha t}\big)^n,
\]
due to the independence assumption. The probability generating function of $Y_k(t)$, assuming that the black ball dies before time $t$, is given 
by 
\begin{equation*}
\begin{split}
&\int_{0}^{t}\Big(\frac{d }{du} (1-e^{ -\delta u}) \Big)\cdot (1 +(v-1) e^{-\alpha(t-u)})^{p} du=\int_{0}^{t}\delta e^{ -\delta u}\cdot (1 +(v-1) e^{-\alpha(t-u)})^p du,
\end{split}
\end{equation*}
due to the fact that $p$ independent white balls are being born at the death of the black ball. Consequently, the probability generating function of the children of the $m-1$ black balls, dying before time $t$, and 
the corresponding number of surviving child balls up to time $t$ is given by 
\begin{equation*}
\begin{split}
&\Big(\int_{0}^{t}\delta e^{ -\delta u}\cdot (1 +(v-1) e^{-\alpha(t-u)})^p du\Big)^{m-1}.
\end{split}
\end{equation*}
Furthermore, the density of the last remaining black ball is given by $\delta\cdot e^{-\delta t}$, giving birth to $p$ more surviving white balls.
Moreover, this final ball can be any one out the $m$ balls. Alltogether, considering all possible final death times $t$, or more precisely 
by conditioning on the stopping time $\tau$ we for obtain the probability generating function $h_{n,m}(v) = \E(v^{X_{n,m}})$ the following result. 

\begin{theorem}
For arbitrary $\alpha,\delta\in\N$ and $p\in \N_0$, the probability generating function of the random variable $X_{n,m}$ counting the number of remaining white balls (divided by $\alpha$) when all black balls have been drawn, $M =\bigl(\begin{smallmatrix} -\alpha & 0 \\ \gamma & -\delta\end{smallmatrix}\bigr)$, $\gamma=\alpha\cdot p$, is given by
\[
h_{n,m}(v) = \int_{0}^{\infty}\big(1+(v-1)e^{-\alpha t}\big)^n \cdot \Big(\int_{0}^{t}\delta e^{ -\delta u}\cdot (1 +(v-1) e^{-\alpha(t-u)})^{p} du\Big)^{m-1} \cdot v^p m\delta   e^{-\delta t}dt.
\]
\end{theorem}
The results above unify and extend the known results of ~\cite{HwaKuPa2007}. 
Moreover, it allows to largely extend the results of~\cite{PanKu2007+} concerning the structure of the moments, as stated below, 
and also to give a complete analysis of the limit laws. Note that by setting $p=0$ one also 
gets the probability generating function for a certain generalized sampling without replacement urn model.
From the result above we will derive a closed formula for the $s$-th factorial moment $\E(\fallfak{\tilde{X}_{n,m}}{s})$ of $\tilde{X}_{n,m}=X_{n,m}-p$, 
for $a\neq d$ such that $a\ell-d\neq 0$; the special case $a=d$ has already been treated in~\cite{HwaKuPa2007}. 
Note that the factorial moments $\E(\fallfak{X_{n,m}}{s})$ of $X_{n,m}$ are recovered using the binomial theorem for the falling factorials
\[
\E(\fallfak{X_{n,m}}{s})=\E\big(\fallfak{(\tilde{X}_{n,m} + p)}{s}\big)=\sum_{\ell=0}^{s}\binom{s}{\ell}\E(\fallfak{\tilde{X}_{n,m}}{\ell})p^{s-\ell}.
\]

\begin{theorem}
The factorial moments of the random variable $\tilde{X}_{n,m}=X_{n,m}-p$ are given in terms of a generalized beta integral,
\begin{equation*}
\begin{split}
\E(\fallfak{\tilde{X}_{n,m}}{s})
&=\delta^{m-1}s!\sum_{j=0}^{s}\binom{n}{j}\sum_{\substack{\sum_{\ell=0}^{p}k_\ell=m-1,\,\sum_{\ell=0}^{p}\ell k_\ell=s-j\\k_\ell\ge 0}}\binom{m-1}{k_0,\dots,k_\ell}\frac{\prod_{\ell=0}^{p}\binom{p}{\ell}^{k_\ell}}{\prod_{\ell=0}^{p}(\ell\alpha-\delta)^{k_\ell}}\\
&\qquad\quad \times\int_{0}^{1}q^{\frac{\alpha j}{\delta}+m-1}\prod_{\ell=0}^{p}(1-q^{\frac{\alpha\ell}{\delta}-1})^{k_\ell}dq.
\end{split}
\end{equation*}
In particular, for we obtain for $p=1$ the simple expression
\[
\E(\fallfak{\tilde{X}_{n,m}}{s})=s!\sum_{\ell=0}^{s}\frac{\binom{n}{\ell}\binom{m}{s-\ell}}{\big(\frac{a}{d}-1\big)^{s-\ell}}
\sum_{i=0}^{s-\ell}(-1)^{s-\ell-i}\frac{\binom{s-\ell}{i}}{\binom{m-s+\ell+i+\frac{a}{d}(s-i)}{m-s+\ell}}.
\]
\end{theorem}

\smallskip

From the result above we will derive a closed formula for the $s$-th factorial moment $\E(\fallfak{\tilde{X}_{n,m}}{s})$ of $\tilde{X}_{n,m}=X_{n,m}-p$, 
for $a\neq d$ such that $a\ell-d\neq 0$; the special case $a=d$ has already been treated in~\cite{HwaKuPa2007}. 
Our starting point is the following expression for $\E(\fallfak{\tilde{X}_{n,m}}{s})$:
\[
\E(\fallfak{\tilde{X}_{n,m}}{s})=E_vD_v^s\frac{h_{n,m}(v)}{v^p},
\]
where $E_v$ denotes the operator which evaluates at $v=1$, and $D_v$ the differentiation operator.
By the binomial theorem we have
\begin{equation*}
\begin{split}
\int_{0}^{t}\delta e^{ -\delta u}\cdot (1 +(v-1) e^{-\alpha(t-u)})^{p} du
&= \delta\sum_{\ell=0}^p\binom{p}{\ell}(v-1)^\ell \frac{e^{ -t\delta}-e^{-\alpha \ell t}}{\ell\alpha-\delta}.  
\end{split}
\end{equation*}
Consequently, using the multinomial theorem, we obtain
\begin{equation*}
\Big(\int_{0}^{t}\delta e^{ -\delta u}\cdot (1 +(v-1) e^{-\alpha(t-u)})^{p} du\Big)^{m-1}
= \delta^{m-1}\sum_{\substack{k_0+\dots+k_p=m-1\\k_\ell\ge 0}}\binom{m-1}{k_0,\dots,k_\ell}
\prod_{\ell=0}^{p}\binom{p}{\ell}^{k_\ell}(v-1)^{\ell k_\ell}.  
\end{equation*}
Using 
\[
\big(1+(v-1)e^{-\alpha t}\big)^n e^{-\delta t}=\sum_{j=0}^{n}\binom{n}{j}(v-1)^je^{-(\alpha j+\delta)t},
\]
we get
\begin{equation*}
\begin{split}
\E(\fallfak{\tilde{X}_{n,m}}{s})
&=\delta^{m}\sum_{j=0}^{n}\binom{n}{j}\sum_{\substack{k_0+\dots+k_p=m-1\\k_\ell\ge 0}}\binom{m-1}{k_0,\dots,k_\ell}\frac{\prod_{\ell=0}^{p}\binom{p}{\ell}^{k_\ell}}{\prod_{\ell=0}^{p}(\ell\alpha-\delta)^{k_\ell}}
E_vD_v^s(v-1)^{j+\sum_{\ell=0}^p \ell k_\ell}\\
&\qquad\quad \times\int_{t=0}^{\infty}e^{-(\alpha j+\delta)t}\prod_{\ell=0}^{p}(e^{ -t\delta}-e^{-\alpha \ell t})^{k_\ell}dt.
\end{split}
\end{equation*}
Since 
\[
E_vD_v^s(v-1)^{j+\sum_{\ell=0}^p \ell k_\ell}=
\begin{cases}
s!, & j+\sum_{\ell=0}^p \ell k_\ell=s,\\
0, & j+\sum_{\ell=0}^p \ell k_\ell\neq s,
\end{cases}
\]
we get the simpler expression
\begin{equation*}
\begin{split}
\E(\fallfak{\tilde{X}_{n,m}}{s})
&=\delta^{m}s!\sum_{j=0}^{s}\binom{n}{j}\sum_{\substack{\sum_{\ell=0}^{p}k_\ell=m-1,\,\sum_{\ell=0}^{p}\ell k_\ell=s-j\\k_\ell\ge 0}}\binom{m-1}{k_0,\dots,k_\ell}\frac{\prod_{\ell=0}^{p}\binom{p}{\ell}^{k_\ell}}{\prod_{\ell=0}^{p}(\ell\alpha-\delta)^{k_\ell}}\\
&\qquad\quad \times\int_{t=0}^{\infty}e^{-(\alpha j+\delta)t}\prod_{\ell=0}^{p}(e^{ -t\delta}-e^{-\alpha \ell t})^{k_\ell}dt.
\end{split}
\end{equation*}
Now we use the substitution $q=e^{-\delta t}$ in order to convert the integral above into a beta-function type integral, 
which proves our result.

\subsection{Higher dimensional urn models}
One can readily extend the $2\times 2$ transition matrix~\eqref{AofA2012-Pills} to higher dimensions,
\begin{small}
\begin{equation*}
    M = \left(
    \begin{smallmatrix}
        -\alpha_1 & 0 & 0 & \cdots & 0 & 0 & 0 \\[-1ex]
        p_2\alpha_1 & -\alpha_2 & 0 & \ddots & \ddots & \ddots & 0 \\[-1ex]
        0 & p_3\alpha_2 & -\alpha_3 & \ddots & \ddots & \ddots & 0 \\[-1ex]
        \vdots & \ddots & \ddots & \ddots & \ddots &
        \ddots & \vdots \\[-1ex]
        0 & \ddots & \ddots & \ddots & -\alpha_{r-2} & 0 & 0 \\[-1ex]
        0 & \ddots & \ddots & \ddots & p_{r-1}\alpha_{r-2} & -\alpha_{r-1} & 0 \\
        0 & 0 & 0 & \cdots & 0 & p_{r}\alpha_{r-1} & -\alpha_r
    \end{smallmatrix}
    \right),
\end{equation*}
\end{small}
with $\alpha_i\in\N$ and $p_i\in\N_0$. We consider the distribution of the random vector $\mathbf{X}_{\mathbf{n}}=(X_{\mathbf{n}}^{[1]},\dots,X_{\mathbf{n}}^{[r-1]})$, which
counts the number of type 1 up to type $r-1$ pills when all pills of $r$ units
are all taken, starting with $n_i$ pills of $i$ units, $i=1,\dots,r$.
One may use similar arguments to the $2\times 2$ case to obtain the following result.

\begin{theorem}
The probability generating function of $\mathbf{X}_{\mathbf{n}}$ is given by
\[
h_{\mathbf{n}}(\mathbf{v}) = \int_{0}^{\infty}\Big(g_{r}(t,\mathbf{v})\Big)^{n_r-1}v_{r-1}^{p_r} \alpha_r n_r e^{-\alpha_r t}\prod_{\ell=1}^{r-1}\big(f_\ell(t,\mathbf{v})\big)^{n_\ell} 
dt.
\]
Here $f_j(t,\mathbf{v})$ denotes a sequence of functions defined by $f_0(t,\mathbf{v})=1$, 
and 
\[
f_j(t,\mathbf{v})=v_je^{-\alpha_j t}+\int_{0}^{t}\alpha_j e^{-\alpha_j u_j} \big(f_{j-1}(t-u_j,\mathbf{v})\big)^{p_j}du_j,\quad j\ge 1,
\]
with $g_j(t,\mathbf{v})=f_j(t,\mathbf{v})-v_je^{-\alpha_j t}$.
\end{theorem}

\subsection{General weight sequences}
One may also obtain the result for $X_{n,m}$ using a slightly different model. Our first process still one consists of $n$ independent ordinary death processes (white balls) with death rate $\alpha$. 
However, concerning the second process, we consider a single modified death process $B(t)$ with death rates $\theta_m,\dots,\theta_1$, starting with $B(0)=m$.
At each transition of $B(t)$ exactly $p$ white balls are being born, modelled by $p$ independent ordinary death processes (white balls) with death rate $\alpha$. 
Consequently, one obtains the alternative description
\[
h_{n,m}(v) = \int_{0}^{\infty}\big(1+(v-1)e^{-\alpha t}\big)^n \cdot p_m(t,v)dt.
\]
where $p_m(t,v)$ denotes density of $B(t)$ dying out before time $t$, with variable $v$ marking the living white balls at time $t$,
\begin{equation}
\begin{split}
p_m(t,v)&=\frac{d}{dt}\bigg(\int_{0}^{t}\theta_m e^{-\theta_m u_m}(1 +(v-1)e^{-\alpha(t-u_m)})^{p}\\
&\quad\qquad\qquad\int_{0}^{t-u_m}\theta_{m-1}  e^{-\theta_{m-1}u_{m-1}}(1 +(v-1)e^{-\alpha(t-u_m-u_{m-1})})^{p}\dots\\
&\qquad\qquad\dots\int_{0}^{t-\sum_{\ell=2}^{m}u_\ell}\theta_{1}  e^{-\theta_1 u_1}(1+(v-1)e^{-\alpha(t-\sum_{\ell=1}^{m}u_\ell)})^{p}du_1\dots du_m\bigg).
\end{split}
\end{equation}
Note that for $\theta_k=\delta\cdot k$, one reobtains the earlier result.

\section{Probabilistic analysis of sampling without replacement and OK Corral type urn models}
We will generalize the sampling without replacement urns, and OK Corral urn models by analyzing two urn models associated to
sequences of positive real numbers $A=(\alpha_n)_{n\in\N}$ and $B=(\beta_m)_{m\in\N}$.
The dynamics of the discrete time process of drawing and replacing balls is as follows:
At every discrete time step, we draw a ball from the urn according to the number of white and black balls present in the urn, with respect to the sequences $(A,B)$, subject to the two models defined below. The choosen ball is discarded and the sampling procedure continues until one type of balls is completely drawn.\\[0.2cm] 
\emph{Urn model I (Sampling with replacement with general weights)}. Assume that $n$ white and $m$ black balls are contained in the urn, with arbitrary $n,m\in \N$. A white ball is drawn with probability $\alpha_n/(\alpha_n+\beta_m)$, and a black ball is drawn with probability $\beta_m/(\alpha_n+\beta_m)$. Additionally, we assume for urn model I that $\alpha_0=\beta_0=0$. \\[0.2cm]
\emph{Urn model II (OK Corral urn model with general weights)}. For arbitrary $n,m\in \N$ assume that $n$ white and $m$ black balls are contained in the urn. A white ball is drawn with probability $\beta_m/(\alpha_n+\beta_m)$, and a black ball is drawn with probability $\alpha_n/(\alpha_n+\beta_m)$. \\[0.2cm]
The absorbing states, i.e.~the points where the evolution of the urn models stop, are given for both urn models by the positive lattice points on the the coordinate axes $\{(0,n)\mid n\ge 1\} \cup \{(m,0)\mid n\ge 1\}$. These two urn models generalize two famous P\'olya-Eggenberger urn models with two types of balls, namely the classical sampling without replacement (I), and the so-called OK-Corral urn model (II), described in detail below. 
We are interested in a probabilisitic derivation of the distribution of the random variable $X_{n,m}$, counting the number of white balls, when all black balls have been drawn. 
In order to simplify the analysis we note that there only exists a single one urn model. 
\begin{lemma}[\cite{Me2010}]
Let $\P\{X_{n,m,[A,B,I]}=k\}$ denote the probability that $k$ white balls remain when all black balls have been drawn in urn model I with weight
sequences $A=(\alpha_n)_{n\in\N}$, $B=(\beta_m)_{m\in\N}$ and $\P\{X_{n,m,[\tilde{A},\tilde{B},II]}=k\}$ the corresponding probability in urn model II with weight
sequences $\tilde{A}=(\tilde{\alpha}_n)_{n\in\N}$, $\tilde{B}=(\tilde{\beta}_m)_{m\in\N}$. The probabilities $\P\{X_{n,m,[A,B,I]}=k\}$ and $\P\{X_{n,m,[\tilde{A},\tilde{B},II]}=k\}$ are dual to each other, i.e.~ they are related in the following way.
\begin{equation*}
\P\{X_{n,m,[A,B,I]}=k\}=\P\{X_{n,m,[\tilde{A},\tilde{B},II]}=k\},
\end{equation*}
for $\alpha_n=\frac{1}{\tilde{\alpha}_n}$, $\beta_m =\frac{1}{\tilde{\beta}_m}$, $n,m\in\N$, and $k>0$.
\end{lemma}
Without loss of generality, we will restrict ourselves to the urn model I. Note that the recurrence relation for the probability generating function 
$h_{n,m}(v)=\E(v^{X_{n,m}})$ of $X_{n,m}$ is given by 
\begin{equation}
\begin{split}
\label{AofA2012-Samp1}
h_{n,m}(v)= \frac{\alpha_n}{\alpha_n+\beta_m}h_{n-1,m}(v) +  \frac{\beta_m}{\alpha_n+\beta_m}h_{n,m-1}(v), \quad n,m\ge 1,
\end{split}
\end{equation}
with initial values $h_{n,0}(v)=v^n$, $h_{0,m}(v)=1$ $n,m\ge 0$.

\subsection{Probabilistic embedding}
We use a probabilistic approach, embedding the discrete-time model into a continuous-time model.
The basic idea is as follows. We consider two independent death processes $X(t)$, and $Y(t)$, which stop at 0. 
Their death rates are are defined using the weight sequences,
$A=(\alpha_n)_{n\in\N}$, $B=(\beta_m)_{m\in\N}$: 
the death rates of $X(t)$, starting with $X(0)=n$ are $\alpha_n,\dots,\alpha_1$, and the death rates of $Y(t)$, starting with $Y(0)=m$ 
are $\beta_m,\dots,\beta_1$. For the sake of convenience we set $\beta_0=0$.
We can model the random variable $X_{n,m}$ of urn model $I$
by looking at the distribution of $C_{n,m}=X(\tau)$, starting with $X(0)=n$, where $\tau$ denotes the time
of the process $Y(t)$ dying out, $\tau=\inf\{t>0: Y(t)=0\}$. 
By conditioning on the first transition of the two processes one directly obtains 
the recurrence relation~\eqref{AofA2012-Samp1} for $\E(v^{C_{n,m}})$, which proves
that $X_{n,m}$ and $C_{n,m}$ have the same distribution. Now things are simple. The probability that the process $X(t)=k$, is according to the definition given by the iterated integral
\begin{equation*}
\begin{split}
\P\{X(t)=k\}&=\int_{0}^{t}\alpha_n e^{-\alpha_n u_n}\int_{0}^{t-u_n}\alpha_{n-1}e^{-\alpha_{n-1} u_{n-1}} \dots \\
&\qquad\qquad \dots\int_{0}^{t-u_n-\dots -u_{k+2}}\alpha_{k+1}e^{-\alpha_{k+1} u_{k+1}}\cdot e^{\alpha_k(t-u_n-\dots -u_{k+2}-u_{k+1})} du_{k+1}\dots du_n.
\end{split}
\end{equation*}
This integral can be evaluated,
\begin{equation}
\begin{split}
\label{AofA2012-OK1}
\P\{X(t)=k\}&= \bigg(\prod_{h=k+1}^{n}\alpha_h\bigg)\sum_{h=k}^{n}\frac{e^{-\alpha_h t}}{\prod_{\substack{j=k\\j\neq h}}^{n} (\alpha_j-\alpha_h)},
\end{split}
\end{equation}
which can easily be checked by induction. This result is covered in standard textbooks or lecture notes, it's derivation is usually based on the Kolomogorov equation
and an application of the Laplace transform. 
The exact distribution of $\tau$ is given by
\[
\P\{\tau< t\}=\int_{0}^{t}\beta_m e^{-\beta_m u_m}du_{m}\int_{0}^{t-u_m}\beta_{m-1}e^{-\beta_{m-1} u_{m-1}} \dots 
\int_{0}^{t-u_m-\dots -u_2}\beta_{1}e^{-\beta_{1} u_{1}}du_1\dots du_m.
\]
One obtains the closed formula
\begin{equation*}
\P\{\tau< t\}=1+\bigg(\prod_{\ell=1}^{m}\beta_{\ell}\bigg)\sum_{\ell=1}^{m}\frac{e^{-\beta_\ell t}}{\prod_{\substack{i=0\\i\neq \ell}}^m (\beta_i-\beta_\ell)},
\end{equation*}
using the convention $\beta_0=0$. Hence, then density function of the stopping time $\tau$ is given by
\begin{equation}
\label{AofA2012-OK2}
\frac{d}{dt}\P\{\tau< t\}=\bigg(\prod_{\ell=1}^{m}\beta_{\ell}\bigg)\sum_{\ell=1}^{m}\frac{e^{-\beta_\ell t}}{\prod_{\substack{i=1\\i\neq \ell}}^m (\beta_i-\beta_\ell)}.
\end{equation}

Considering all possible times when the second process dies out 
leads to the integral representation
\begin{equation}
\begin{split}
\label{AofA2012-Samp2}
\P\{X_{n,m}=k\}&=\int_{0}^{\infty}\frac{d}{dt}\P\{\tau< t\}\P\{X(t)=k\}dt \\
&=\bigg(\prod_{\ell=1}^{m}\beta_{\ell}\bigg)\bigg(\prod_{h=k+1}^{n}\alpha_h\bigg)\sum_{h=k}^{n}\sum_{\ell=1}^{m}
\frac{1}{\prod_{\substack{i=1\\i\neq \ell}}^m (\beta_i-\beta_\ell)\prod_{\substack{j=k\\j\neq h}}^{n} (\alpha_j-\alpha_h)}\displaystyle{\int_{0}^{\infty}e^{-(\beta_\ell +\alpha_h) t}dt}\\
&=\bigg(\prod_{\ell=1}^{m}\beta_{\ell}\bigg)\bigg(\prod_{h=k+1}^{n}\alpha_h\bigg)\sum_{h=k}^{n}\sum_{\ell=1}^{m}
\frac{1}{(\beta_\ell +\alpha_h)\prod_{\substack{i=1\\i\neq \ell}}^m (\beta_i-\beta_\ell)\prod_{\substack{j=k\\j\neq h}}^{n} (\alpha_j-\alpha_h)}.
\end{split}
\end{equation}
The result above can be simplified in two different ways using the partial fraction identities
\begin{equation}
\begin{split}
\label{AofA2012parfrac}
\frac{1}{\prod_{j=k}^{n}(\alpha_j+x)}&= \sum_{h=k}^n\frac{1}{(x+\alpha_h)\prod_{\substack{j=k\\j\neq h}}^{n}(\alpha_j-\alpha_h)},\\
\frac{1}{\prod_{i=1}^{m}(\beta_i+x)}&= \sum_{\ell=1}^n\frac{1}{(x+\beta_\ell)\prod_{\substack{i=1\\i\neq \ell}}^{m}(\beta_i-\beta_\ell)}.
\end{split}
\end{equation}
Consequently, we have obtain a transparent probabilistic proof of the following result.
\begin{theorem}[\cite{Me2010}]
The probability mass function of the random variable $X_{n,m}$, counting the number of remaining white balls when all black balls
have been drawn in urn model I with weight sequences $A=(\alpha_n)_{n\in\N}$, $B=(\beta_m)_{m\in\N}$, is for $n,m\ge 1$ and $n\ge k\ge 1$ given by the explicit formula 
\begin{equation*}
\begin{split}
\P\{X_{n,m}=k\}&= \Big(\prod_{h=1}^{m}\beta_h\Big)\Big(\prod_{h=k+1}^{n}\alpha_h\Big)\sum_{\ell=1}^{m}
\frac{1}{\Big(\prod_{j=k}^{n}(\alpha_j+\beta_{\ell})\Big)\Big(\prod_{\substack{i=1\\i\neq \ell}}^{m}(\beta_{i}-\beta_{\ell})\Big)}\\
&= \Big(\prod_{h=1}^{m}\beta_h\Big)\Big(\prod_{h=k+1}^{n}\alpha_h\Big)\sum_{\ell=k}^{n}
\frac{1}{\Big(\prod_{\substack{j=k\\j\neq \ell}}^{n}(\alpha_j-\alpha_{\ell})\Big)\Big(\prod_{i=1}^{m}(\beta_{i}+\alpha_{\ell})\Big)},
\end{split}
\end{equation*}
assuming that $\alpha_j \neq \alpha_{\ell}$ and $\beta_j \neq \beta_{\ell}$, $1\le j<\ell<\infty$, and that $\alpha_0=0$.
\end{theorem}
It can be shown that the result above is also valid for $k=0$. Moreover, by the duality of the two urn models, one also gets the corresponding result for the urn model II, OK-Corral type urn models, by switching to weight sequences,
$\tilde{A}=(1/\alpha_n)_{n\in\N}$, $\tilde{B}=(1/\beta_m)_{m\in\N}$. 

\subsection{Sums of independent exponential random variables}
Of course, the formulas~\eqref{AofA2012-OK1},~\eqref{AofA2012-OK2} stated before do not come as a surprise, since one can take yet another viewpoint. 
The time $\tau$ until the second process $Y(t)$ dies out has the same distribution as the sum of $m$ independent exponential distributed random variables $\epsilon_{\beta_i}$
with parameters $\beta_m,\dots,\beta_1$ stemming from the death rates of the process. 
Hence, 
\[
\tau=\bigoplus_{\ell=1}^{m}\epsilon_{\beta_\ell},
\] 
where $\epsilon_{\beta_i}$ denotes an exponential distribution with parameter $\beta_i$,
and the density is simply the formula stated in~\eqref{AofA2012-OK2}.
Furthermore, the distribution of $X(t)$ can also be modelled by $k$ independent random variables:
let 
\[
\theta=\bigoplus_{\ell=k+1}^{n}\epsilon_{\alpha_\ell}.
\]
If $\P\{X(t)=k\}$, then the $k$ transitions of the process $X$ have occured before $t$, 
and no more transition afterwards. Hence, 
\begin{equation*}
\begin{split}
\P\{X(t)=k\}&=\P\{\theta<t, \theta+\epsilon_{\alpha_k}>t\}
=\int_{0}^{t}\bigg(\prod_{h=k+1}^{n}\alpha_h\bigg)\sum_{h=k+1}^{n}\frac{e^{-\alpha_h u}}{\prod_{\substack{j=k+1\\j\neq h}}^{n} (\alpha_j-\alpha_h)}e^{-\alpha_k(t-u)}du\\
&=\sum_{h=k+1}^{n}\frac{\bigg(\prod_{h=k+1}^{n}\alpha_h\bigg)e^{-\alpha_h t}}{\prod_{\substack{j=k\\j\neq h}}^{n} (\alpha_j-\alpha_h)}
+\sum_{h=k+1}^{n}\frac{\bigg(\prod_{h=k+1}^{n}\alpha_h\bigg)e^{-\alpha_k t}}{(\alpha_h-\alpha_k) \prod_{\substack{j=k+1\\j\neq h}}^{n} (\alpha_j-\alpha_h)}
\end{split}
\end{equation*}
which simplifies to~\eqref{AofA2012-OK1} after an application of~\eqref{AofA2012parfrac}.

\end{document}